\documentclass[12pt,a4paper]{article}

\title{An Indefinite Convection-Diffusion Operator With Real Spectrum}
\author{John Weir}

\newtheorem{theorem}{Theorem}[section]
\newtheorem{lemma}[theorem]{Lemma}
\newtheorem{proposition}[theorem]{Proposition}

\newenvironment{proof}[1][Proof]{\begin{trivlist}
\item[\hskip \labelsep {\bfseries #1}]}{\end{trivlist}}

\newcommand{\qed}{\nobreak \ifvmode \relax \else
      \ifdim\lastskip<1.5em \hskip-\lastskip
      \hskip1.5em plus0em minus0.5em \fi \nobreak
      \vrule height0.75em width0.5em depth0.25em\fi}

\newcommand{\setR}{\mathbf{R}}
\newcommand{\setZ}{\mathbf{Z}}
\newcommand{\setC}{\mathbf{C}}
\newcommand{\eqnref}[1]{(\ref{#1})}
\newcommand{\norm}[1]{\left|\left|#1\right|\right|}
\newcommand{\dom}[1]{\mathrm{Dom}\left(#1\right)}

\newcommand{\e}{\mathrm{e}}

\begin{document}
\maketitle

\section{Introduction}

For $0 < \varepsilon < 2$, we consider the operator
\begin{equation}
(Hf)(\theta) := \varepsilon \frac{\partial}{\partial \theta} \left ( \sin(\theta) \frac{\partial f}{\partial \theta} \right ) + \frac{\partial f}{\partial \theta}
\end{equation}
initially defined on all $\mathcal{C}^2$ periodic functions $f \in L^2(-\pi, \pi)$; the exact domain is given by taking the closure of the operator defined on the above functions. That such a closure exists will be shown later.
In a recent paper \cite{bobs} Benilov, O'Brien and Sazonov showed that the equation
\begin{equation}
\frac{\partial f}{\partial t} = Hf
\end{equation}
approximates the evolution of a liquid film inside a rotating horizontal cylinder.

We shall show that the eigenvalue problem
\begin{equation}\label{eqn:eig}
-iHf = \lambda f
\end{equation}
has only real eigenvalues, which were shown to exist by Davies in \cite{davies-2007}. In the same paper, he showed that the spectrum of $-iH$ is equal to the set of its eigenvalues, so this implies that the spectrum is real. This was conjectured in \cite{bobs}, and Chugunova and Pelinovsky proved in \cite{chugunova-2007} that all but finitely many eigenvalues are real, and gave numerical evidence that all eigenvalues are real. Our approach is to show that it is sufficient to consider $H$ acting on the Hardy space $H^2(-\pi,\pi)$ and analytically continue any solution of \eqnref{eqn:eig} to the unit disc, where the corresponding ODE \eqnref{eq:sl} is now self-adjoint on $[0,1]$ with regular singularities at the end points. In order to do this we make use of a bound on the Fourier coefficients proved by Davies in \cite{davies-2007}.

As in \cite{davies-2007}, by expanding $f \in L^2(- \pi, \pi)$ in the form
\[ f(\theta) = \frac{1}{\sqrt{2 \pi}} \sum_{n \in \setZ} v_n \e^{in\theta}, \]
one may rewrite the eigenvalue problem in the form $Av = \lambda v$, where $A = -iH$ is given by
\[ (Av)_n = \frac{\varepsilon}{2} n (n-1) v_{n-1} - \frac{\varepsilon}{2} n (n+1) v_{n+1} + n v_n. \]
Here we have identified $l^2(\setZ)$ and $L^2(-\pi,\pi)$ using the Fourier transform $\mathcal{F} : l^2(\setZ) \to L^2 (-\pi, \pi)$. We have
\[A = -i\mathcal{F}^{-1}H\mathcal{F}\]
and define $\dom{H} = \mathcal{F}(\dom{A})$.

The (unbounded) tridiagonal matrix $A$ is of the form
\[ A = \left (	\begin{array}{ccc} A_-	& 0 & 0\\
																	0		&	0	&	0\\
																	0		&	0	&	A_+
								\end{array} \right )\]
where $A_-$ acts in $l^2(\setZ_-)$, the central $0$ acts in $\setC$ and $A_+$ acts in $l^2(\setZ_+)$.
We assume that $A_+$ has its natural maximal domain
\[ \mathcal{D} = \{v \in l^2(\setZ_+):A_+ v \in l^2(\setZ_+) \}. \]
Davies has shown in \cite{davies-2007} that $A_+$ is closed and that $\mathcal{D}$ is the closure under the graph norm of $A_+$ of the subspace consisting of those $v \in l^2(\setZ_+)$ that have finite support. Let $\tau$ be the natural identification between $l^2(\setZ_+)$ and $l^2(\setZ_-)$; then $\dom{A} = \tau(\mathcal{D}) \oplus \setC \oplus \mathcal{D}$ and $\tau$ induces a unitary equivalence between $A_+$ and $A_-$. Therefore, in order to prove that all eigenvalues of $A$ are real, we only need to prove that all eigenvalues of $A_+$ are real. The Fourier transform identifies $l^2(\setZ_+)$ with $\{f \in H^2(-\pi,\pi) : \int_{-\pi}^{\pi}f(\theta)\mathrm{d}\theta = 0\}$.

Let $H_0$ be the restriction of $H$ to $\mathcal{C}^2_{\mathrm{per}}([-\pi,\pi])$, which is clearly a subspace of $\dom{H}$. We now show that $H$ is the closure of $H_0$.

\begin{proposition}
Where $H$ and $H_0$ are as above, $H$ is the closure of $H_0$.
\end{proposition}
\begin{proof} 
It follows from Davies' result on the domain of $A_+$ that the trigonometric polynomials are dense in $\dom{H}$ with respect to the graph norm. Since the trigonometric polynomials are contained in $\mathcal{C}^2_{\mathrm{per}}([-\pi,\pi])$, this space is also dense in $\dom{H}$, which is closed in graph norm since $\dom{A}$ is.
\qed
\end{proof}

\section{Reality Of The Eigenvalues}

If $A_+v = \lambda v$, then $v$ is a solution of the recurrence relation
\begin{equation}\label{eqn:recurrence}
\frac{\varepsilon}{2}(n+1)(n+2)v_{n+2} + (\lambda - n-1)v_{n+1} - \frac{\varepsilon}{2}n(n+1)v_n = 0
\end{equation}
satisfying the initial condition $\varepsilon v_2 = (1-\lambda)v_1$. We shall study the generating function, $\sum_{k=1}^{\infty}v_k z^k$ of $(v_k)$.

\begin{lemma}\label{genfn}Let $v \in l^2(\setZ_+)$ be such that $A_+v = \lambda v$. Then the function $u(z) := \sum_{k=1}^{\infty}v_k z^k$, defined for $|z| < 1$, satisfies the differential equation
\begin{equation}
	u'' - 2 \frac{z+1/\varepsilon}{(1-z)(1+z)}u' + \frac{2\lambda/\varepsilon}{z(1-z)(1+z)}u = 0.
	\label{eq:genfn}
\end{equation}
\end{lemma}
\begin{proof}
The constant term in
\[ z(1-z)(1+z)u'' - 2(z+1/\varepsilon)zu' + (2\lambda/\varepsilon)u\]
is clearly $0$, and the coefficient of $z$ is
\[ 2v_2 - 2v_1/\varepsilon + 2\lambda v_1/\varepsilon = 2 \left(v_2 - \frac{(1-\lambda)}{\varepsilon}v_1\right) = 0.\]
The coefficient of $z^n$ is
\begin{eqnarray*}
	n(n+1)v_{n+1} - (n-1)(n-2)v_{n-1} - \frac{2}{\varepsilon}nv_n-2(n-1)v_{n-1} + \frac{2\lambda}{\varepsilon}v_n \\
	=\frac{2}{\varepsilon}\left[\frac{\varepsilon}{2}n(n+1)v_{n+1}+(\lambda-n)v_n -\frac{\varepsilon}{2}n(n-1)v_{n-1}\right] = 0
\end{eqnarray*}
for $n \geq 2$.
\qed
\end{proof}

From here on we assume that $\lambda \in \setC$ is an eigenvalue of the operator $A_+$ and $v$ is a corresponding non-zero eigenvector. In \cite{davies-2007}, Davies proved that there exist constants $b, m$ such that
\begin{equation}
	\norm{v}_{\infty,c} \leq b |\lambda|^m \norm{v}_2,
	\label{eq:daviesbound}
\end{equation}
where $c = 1+1/\varepsilon$ and $\norm{w}_{\infty,c} := \sup\{|w_n| n^{c} : 1 \leq n < \infty \}$. Since $\varepsilon > 0$, this implies that $v \in l^1(\setZ_+)$, and hence that $\sum_{k=1}^{\infty}v_k z^k$ is absolutely convergent for $|z| \leq 1$. The equation \eqnref{eq:genfn} can be written in the form of Heun's equation
\begin{equation}
	u'' + \left(\frac{\gamma}{z} + \frac{\delta}{z-1} + \frac{\epsilon}{z-a}\right)u' + \frac{\alpha\beta z - \mu}{z(z-1)(z-a)}u = 0
	\label{eq:heun}
\end{equation}
with $\alpha = 1$, $\beta=0$, $\gamma = 0$, $\delta = 1+1/\varepsilon$, $\epsilon = 1-1/\varepsilon$, $a=-1$ and $\mu = 2\lambda/\varepsilon$. This is a Fuchsian equation with four regular singular points, at $0,1,-1$ and $\infty$, with $\{0,1-\gamma\}$, $\{0,1-\delta\}$, $\{0,1-\epsilon\}$ and $\{\alpha, \beta\}$ as the roots of the corresponding indicial equations (for the Frobenius series at each regular singular point). For background information on Heun's equation, see \cite{ronveaux}.

\begin{lemma}\label{lem:heun}Suppose that $0< \varepsilon <2$, $1/ \varepsilon \notin \setZ$, $\lambda \in \setC$ is an eigenvalue of $A_+$ and $v \in l^2(\setZ_+)$ is a corresponding non-zero eigenvector. Then there exists a solution $u$ of \eqnref{eq:heun} which is analytic in an open set containing $[0,1]$ and such that $u(z) = \sum_{k=1}^{\infty}v_k z^k$ for all $z$ such that $|z| \leq 1$.
\end{lemma}
\begin{proof}
Define $u(z) = \sum_{k=1}^{\infty}v_k z^k$ on $\{z \in \setC : |z| < 1\}$. Let $u_1$ be the solution of \eqnref{eq:heun} with exponent $0$ about $1$ and $u_2$ be the solution with exponent $-1/\varepsilon$ about $1$. Let $U$ be the intersection of the open discs of unit radius about $0$ and $1$. The space of solutions of \eqnref{eq:heun} in $U$ is two-dimensional, and $u$, $u_1$ and $u_2$ lie in this space. Hence there exist constants $a$, $b$ such that $u = a u_1 + b u_2$ in $U$. Since $v \in l^1(\setZ_+)$, $u(z)$ converges to a finite limit as $z \to 1$ in $U$. Also $u_1(1)$ is finite, but $u_2(z) \to \infty$ as $z \to 1$ in $U$. Therefore we must have $b=0$ and $u = a u_1$ in $U$. Let $W$ be the union of the open discs of unit radius about $0$ and $1$. We now extend $u$ to all of $W$ by $u(z) = a u_1(z)$ on the open disc of unit radius about $1$. Now $u$ is an analytic solution of \eqnref{eq:heun} on $W$, which is an open set containing $[0,1]$, such that $u(z) = \sum_{k=1}^{\infty}v_k z^k$ for all $z$ such that $|z| \leq 1$.
\qed
\end{proof}

\begin{theorem}\label{thm:realspectrum} Suppose that $0< \varepsilon <2$, $1/ \varepsilon \notin \setZ$ and $\lambda \in \setC$ is an eigenvalue of $A_+$. Then $\lambda \in \setR$.
\end{theorem}
\begin{proof}
Let $v \in l^2(\setZ_+)$ be a non-zero eigenvector corresponding to $\lambda$. Let $u$ be as in Lemma \ref{lem:heun} and put $\mu = 2 \lambda/\varepsilon$. Following \cite{everitt}, the equation \eqnref{eq:heun} can be written as
\begin{equation}\label{eq:sl}
-(pu')' + qu = \mu w u
\end{equation}
on the complex plane cut along $[1,\infty)$ and $(-\infty,0]$, where
\[\begin{array}{lcccl}
	p(z)	& = & z^{\gamma} (1-z)^{\delta} (z-a)^{\epsilon}									& = & (1-z)^{1+1/\varepsilon}(z+1)^{1-1/\varepsilon}\\
	q(z)	& = & \alpha \beta z^{\gamma} (1-z)^{\delta-1} (z-a)^{\epsilon-1}	& = & 0\\
	w(z)	& = & z^{\gamma-1} (1-z)^{\delta-1} (z-a)^{\epsilon-1}						& = & z^{-1} (1-z)^{1/\varepsilon} (z+1)^{-1/\varepsilon}.
\end{array}\]
We now restrict $u$ to $[0,1]$. It is clear that $u \in \mathcal{C}^{\infty}([0,1])$ and
\[ -(pu')' = \mu w u \]
on $(0,1)$. Note that $p>0$ on $[0,1)$ and $w>0$ on $(0,1)$ with $p(x),w(x) \to 0$ as $x \to 1$ from below. Since $u$ has a zero of order $1$ at $0$ and $w$ has a pole of order $1$ at $0$, $wu \in \mathcal{C}([0,1])$. Therefore $|u|^2w \in L^1(0,1)$. Now
\begin{eqnarray*}
	\mu \int_0^1|u(x)|^2w(x)\mathrm{d}x & = &
	\lim_{n \to \infty}\left\{-\int_{1/n}^{1-1/n}(pu')'(x)\overline{u(x)}\mathrm{d}x\right\}\\
	& = & \lim_{n \to \infty}\left\{-\left[p(x)u'(x)\overline{u(x)}\right]_{1/n}^{1-1/n} + \int_{1/n}^{1-1/n}p(x)u'(x)\overline{u'(x)}\mathrm{d}x\right\}\\
	& = & \lim_{n \to \infty}\left\{\left[p(x)\left\{u(x)\overline{u'(x)}-u'(x)\overline{u(x)}\right\}\right]_{1/n}^{1-1/n} - \int_{1/n}^{1-1/n}\overline{(pu')'(x)}u(x)\mathrm{d}x\right\}\\
	& = & -\int_0^1 \overline{(pu')'(x)}u(x)\mathrm{d}x\\
	& = & \overline{\mu}\int_0^1 |u(x)|^2 w(x)\mathrm{d}x.
\end{eqnarray*}
Since $u$ is a non-zero solution of \eqnref{eq:sl} and $w>0$ a.e. we have $\mu = \overline{\mu}$ and hence $\mu \in \setR$. Since $\lambda = \frac{\varepsilon}{2}\mu$, we also have $\lambda \in \setR$. \qed
\end{proof}

\bibliographystyle{plain}
\bibliography{References}

\end{document}